\documentclass[a4paper,10pt,psamsfonts,reqno]{amsart}
\usepackage{amsmath, amssymb, amscd, amsthm, euscript, verbatim}
\usepackage{epsfig}
\pdfoutput=1

\newtheoremstyle{thmbold}{11.5pt}{}{\slshape}{}{\bfseries}{.}{.39em}{}
\theoremstyle{thmbold}

\newtheorem{thm}{Theorem}
\newtheorem{cor}[thm]{Corollary}
\newtheorem{prop}[thm]{Proposition}
\newtheorem{lem}[thm]{Lemma}

\newtheoremstyle{excap}{\topsep}{10pt}{}{}{\scshape}{.}{.7em}{}
\theoremstyle{excap}

\newcommand{\abs}[1]{\lvert #1 \rvert}
\newcommand{\floor}[1]{\left\lfloor #1 \right\rfloor}

\newcommand{\set}[1]{\{ #1 \}}
\newcommand{\setpres}[2]{\left\{ #1 \mid \text{#2} \right\}}

\newcommand{\card}[1]{\abs{#1}}

\newcommand{\Dbb}{\mathbb{D}}

\newcommand{\Nbb}{\mathbb{N}}

\newcommand{\Zbb}{\mathbb{Z}}

\newcommand{\Acal}{{\mathcal A}}
\newcommand{\Bcal}{{\mathcal B}}
\newcommand{\Gcal}{{\mathcal G}}
\newcommand{\Pcal}{{\mathcal P}}
\newcommand{\Scal}{{\mathcal S}}

\newcommand{\Asymm}{{\mathcal A}^{\,\sigma}}
\newcommand{\Bsymm}{{\mathcal B}^{\,\sigma}}
\newcommand{\Gsymm}{{\mathcal G}^{\,\sigma}}
\newcommand{\Psymm}{{\mathcal P}^{\,\sigma}}
\newcommand{\Ssymm}{{\mathcal S}^{\,\sigma}}
\newcommand{\Asigma}{{A}^{\,\sigma}}
\newcommand{\Asp}[1]{{A^{\,\sigma\,\prime}_{#1}}}

\newcommand{\bsymm}[1]{\beta^{\,\sigma}_{#1}}
\newcommand{\gsymm}[1]{\gamma^{\,\sigma}_{#1}}

\begin{document}


\title[Counting numerical sets with no small atoms]{Counting numerical sets\\ with no small atoms}

\author{Jeremy Marzuola}
\address{Jeremy Marzuola\\ Department of Applied Physics and Applied Mathematics\\ \break Columbia University \\
New York\\ NY \\ 10027} 
\email{jm3058@columbia.edu}
\author{Andy Miller}
\address{Andy Miller\\ Department of Mathematics\\ University of Oklahoma\\ Norman\\ OK\\ 73019}
\email{amiller@math.ou.edu}


\begin{abstract}
A numerical set $S$ with Frobenius number $g$ is a set of integers with $\min(S) = 0$ 
and $\max(\Zbb - S)=g$, and its atom monoid is 
$A(S) = \setpres{n \in \Zbb}{$n+s \in S$ for all $s \in S$}$.
Let $\gamma_g$ be the number of numerical sets $S$ having $A(S) = \set{0} \cup (g,\infty)$
divided by the total number of numerical sets with Frobenius number $g$. We show that
the sequence $\set{\gamma_g}$ is decreasing and converges to a number
$\gamma_\infty \approx .4844$ (with accuracy to within $.0050$). 
We also examine the singularities of the generating function
for $\set{\gamma_g}$. Parallel results are obtained for the ratio 
$\gsymm{g}$ of the number of symmetric numerical sets $S$ with $A(S) = \set{0} \cup (g,\infty)$
by the number of symmetric numerical sets with Frobenius number $g$.
These results yield
information regarding the asymptotic behavior of the number of finite
additive 2-bases.
\end{abstract}

\maketitle

Let $\Zbb$ denote the additive group of integers and let $\Nbb$ 
denote the monoid of nonnegative integers.
Both of these sets are linearly ordered by the Archimedean ordering and we will use 
standard interval notation to describe their convex subsets.
If $n \in \Zbb$ and $S \subseteq \Zbb$ then the translate of $S$ by $n$ is the set
$n+S = \setpres{n+s}{$s \in S$}$.

A {\it numerical set $S$} is  a cofinite subset of $\Nbb$ which contains $0$, and 
its {\it Frobenius number} is the maximal element in the complement $\Nbb - S$.\footnote{This 
definition differs from that employed in \cite{AM} where a `numerical set' would be a translate 
$n+S$ of a numerical set $S$ (in the sense given here) by an arbitrary integer $n$. 
Since the atom monoid of $n+S$ equals the atom monoid of $S$, this variation of the definition 
should not lead to any confusion.} 
Equivalently, a numerical set $S$ with Frobenius number $g$ is a set of integers with $\min(S) = 0$ 
and $\max(\Zbb - S)=g$.
A numerical set which is closed under addition is called a {\it numerical monoid}. 
Every numerical set $S$ has an associated {\it atom monoid $A(S)$} defined by
  $$ A(S) = \setpres{n\in\Zbb}{$n+S \subseteq S$} ,$$
and this is easily seen to be a numerical monoid with the same Frobenius number as $S$. 
Note that $A(S) \subseteq S$ and that $S$ is a numerical monoid if and only if $A(S) =S$.
The nonzero elements of $A(S)$ are referred to as the {\it atoms of $S$}. 

For each $g \ge 0$ let $\Nbb_g$ be the numerical monoid
  $$\Nbb_g = \Nbb - [1,g] = \set{0} \cup (g,\infty) \ ,$$
which has Frobenius number $g$ when $g>0$.\footnote{The 
Frobenius number of $\Nbb_0=\Nbb$ is $-1$, and this is the only numerical set with nonpositive 
Frobenius number.} 
The atom monoid of every numerical set $S$ with Frobenius number $g$ contains
$\Nbb_g$ and the complement $S - \Nbb_g$ is a subset of $(0,g)$.
Conversely, the union of $\Nbb_g$ with any subset of $(0,g)$ is a numerical set with Frobenius number
$g$. Therefore the set
  $$\Scal(g) = \setpres{S \subseteq \Nbb}{$S$ is a numerical set with Frobenius number $g$}$$
is in one-to-one correspondence with the power set $\Pcal(0,g)$ consisting of all subsets of $(0,g)$, and 
$\Scal(g)$ has cardinality $2^{g-1}$. The subset of $\Scal(g)$ consisting of numerical monoids is a much 
more difficult set to enumerate. This is examined in Backelin's paper \cite{B} 
where it is shown that for large values of $g$ roughly 
$3 \times 2^{\floor{(g-1)/2}}$ of the $2^{g-1}$ elements of 
$\Scal(g)$ are numerical monoids.

If $M \in \Scal(g)$ is a numerical monoid then the
{\it anti-atom set of $M$}  
is the set
  $$\Gcal(M) \ = \ \setpres{S\in \Scal(g)}{$A(S)=M$} .$$
This is contained in the larger set $\Scal(M) = \setpres{S\in \Scal(g)}{$M\subseteq A(S)$}$ 
whose elements might be considered to be `$M$-modules'.\footnote{In \cite{BF} the elements
of $\Scal(M)$ are called `relative ideals over $M$'.} Notice that $\Scal(g) = \Scal(\Nbb_g)$ 
and we will also write
$\Gcal(g) = \Gcal(\Nbb_g)$. This paper is motivated by the following question 
which we shall refer to as the {\it Anti-Atom Problem}. 
\begin{center}\label{test}
  {\it For a  given numerical monoid $M$ with Frobenius number $g$\\ how many numerical sets 
    in $\Scal(g)$ have atom monoid $M$?}
\end{center}
Thus, for a given monoid $M$, the Anti-Atom Problem asks to compute the cardinality of
$\Gcal(M)$. 
This problem is certainly unwieldy
given that it fundamentally presupposes an enumeration of the set of numerical monoids in
$\Scal(g)$---an enumeration which Backelin has shown to be intractable at best.
Nevertheless we will be able to frame aspects of the problem
in a clearer light. Our main result will show 
that there is one monoid $M$ in $\Scal(g)$ (that monoid being $M=\Nbb_g$) 
which itself is the atom monoid for approximately  $48.4\%$ of all numerical 
sets in $\Scal(g)$  for large values of $g$.
In order to describe
this in more depth we first need to discuss symmetry and pseudosymmetry
in numerical sets. These concepts are important throughout much of the theory 
of numerical monoids and numerical sets (see \cite{FGH}, \cite{AM} and \cite{A} for example), and 
will play a role in many of our discussions.

A numerical set $S \in \Scal(g)$ is {\it symmetric} if an integer $x$ is an element of
$S$ if and  only if $g-x$ is not an element of $S$. In other words, 
$S$ is symmetric when the reflection on 
$\Zbb$ given by $x \mapsto g-x$ carries $S$ onto its complement $\Zbb - S$. 
Notice that only numerical semigroups with odd Frobenius number can be symmetric.
A numerical set with even Frobenius number $g$ is said to be {\it pseudosymmetric} if 
$g/2 \notin S$ and for each integer $x\ne g/2$, $x$ is an element of $S$ if and only if 
$g-x$ is not an element of $S$. 
Symmetry and pseudosymmetry can also be described using the notion of duality
of numerical sets. If $S \in \Scal(g)$ then the {\it dual of $S$} is the numerical set
$S^* = \setpres{n \in \Zbb}{$g-n \notin S$}$, and it is not hard to show
that $S^* \in \Scal(g)$ and that $A(S^*)=A(S)$ (more background can be found in 
section~1 of \cite{AM}).
The numerical set $S$ is symmetric if and only if $S^* = S$, and it
is pseudosymmetric if and only if $g$ is even and $S^* = S \cup \set{g/2}$.\footnote{
More generally, if the symmetric difference of $S$ and $S^*$ contains no more
than one element then $S$ is symmetric, pseudosymmetric or ``dually pseudosymmetric''
(meaning that $S^*$ is pseudosymmetric).} 
For each numerical set $S \in \Scal(g)$ there is a rational number $type(S)$
no smaller than one,
called the `type of $S$', which satisfies the property that $S$ is symmetric 
if and only if $type(S) = 1$. The type of a numerical monoid $M \in \Scal(g)$ is 
always an integer, and it equals the cardinality of its {\it omitted atom set} 
${\mathcal O}(M) = \setpres{n \in \Zbb - M}{$n + \big( M- \set{0}\big) \subseteq M$}$.
Since ${\mathcal O}(M) \subset \Nbb$ and
$g \in {\mathcal O}(M)$, the type of 
a numerical monoid $M \in \Scal(g)$ is in the interval $[1,g]$, 
and the largest possible value $type(M) = g$ is only achieved when $M=\Nbb_g$.
The following elementary results allow us to solve the Anti-Atom Problem for symmetric
and pseudosymmetric numerical monoids.

\begin{prop} Suppose that $M$ is a numerical monoid and that $S$ is a numerical set 
with $A(S)=M$. Then $M \subseteq S \subseteq M^*$.
\end{prop}

\begin{proof}
Let $S$ be a numerical set in $\Scal(g)$ with $A(S)=M$  and $s \in S$. 
If $g-s$ were an element of $M$ then $g = s + (g-s)$ would be an element of $S$,
which contradicts $g$ being the Frobenius number of $S$. 
Thus $g-s \notin M$ which implies that $s \in M^*$,  and $M = A(S) \subseteq S \subseteq M^*$. 
\end{proof}

\begin{cor} A numerical monoid $M \in \Scal(g)$ is symmetric if and only if 
there is just one numerical set (which must be $M$ itself) whose atom monoid is $M$.
If $M$ is a pseudosymmetric numerical monoid then there are precisely two numerical
sets (which must be $M$ and $M^*$) whose atom monoid is $M$.
\end{cor}

\begin{proof} Let $M \in \Scal(g)$ be a monoid.
If $M$ is not symmetric then $M \ne M^*$ but
$A(M^*)=A(M)=M$, and so there are at least two distinct numerical sets  
in $\Gcal(M)$. On the other hand, if $M$ is symmetric and $S \in \Gcal(M)$
then $M \subseteq S \subseteq M^* = M$ and $S=M$. If $M$ is pseudosymmetric
and $S \in \Gcal(M)$
then $M \subseteq S \subseteq M^* = M \cup \set{g/2}$, so that $S$ equals $M$ or $M^*$.
\end{proof}

This corollary then provides the first positive answers to the Anti-Atom Problem: namely, that 
$\card{\Gcal(M)} = 1$ when $M$ is symmetric and that $\card{\Gcal(M)} = 2$ when $M$ is 
pseudosymmetric.\footnote{We showed that $\card{\Gcal(M)} = 1$ if and only if $M$ is symmetric,
but it is not hard to construct numerical monoids $M$ with
$\card{\Gcal(M)} = 2$ that are not pseudosymmetric.}  
At the other end of the spectrum, we shall show that the 
anti-atom set of $\Nbb_g$ (which is the numerical monoid in $\Scal(g)$ farthest 
removed from being symmetric since $type(\Nbb_g) = g$ is the largest possible
type among all monoids in $\Scal(g)$)
is an order of magnitude larger in size than that of
any other numerical monoid with Frobenius number $g$.
To establish this we will examine the sequence $\gamma_g = \card{\Gcal(g)}/\card{\Scal(g)}$.
We introduce a combinatorially defined sequence of positive integers $\set{A_k}$  with the property that 
$1-\gamma_g$ is a partial sum of the convergent infinite series $\sum_{k=1}^\infty \, A_k 4^{-k}$.
This allows us to show that $\set{\gamma_g}$ is a decreasing convergent sequence and 
that the limit $\gamma_\infty$ is approximately equal to $.484451$, give or take $.0050$.
The integers ${A_k}$ are combinatorially related to integers
${A_k^\prime}$ which turn out to equal $\card{\Gcal(k)}$ (theorem~\ref{thm:GoodSet})
and there is a nice recursive relation between these two sequences (theorem~\ref{thm:recursions}).
This relationship  will enable us to obtain information about the singularities of the
generating functions for the sequences $\set{A_k}$ and $\set{A_k^\prime}$. 

In addition to forming a large subset of $\Scal(g)$, the numerical sets in 
$\Gcal(g)$ have nice properties in terms of the direct sum decompositions discussed in \cite{AM}.
Given numerical sets $S$ and $T$ and relatively prime atoms $a\in A(S)$ and $b \in A(T)$ 
the {\it direct sum of $S$ and $T$} is the numerical set $bS \oplus aT= \setpres{bs+at}{$s \in S$ 
and $t\in T$}$. Every numerical set $S$ can be trivially described as $S = 1 S \oplus a \Nbb$ 
for any nonzero $a \in A(S)$, but if this is the only kind of direct sum decomposition of $S$ 
then we say that $S$ is {\it irreducible}. Every numerical set can be expressed as a finite 
direct sum of irreducibles. By \cite[Proposition~4.4]{AM}, the only numerical set in ${
\bigcup \setpres{{\Gcal}(g)}{$g \ge 1$}}$
which is not irreducible is $\Nbb_1 = 2\Nbb \oplus 3\Nbb$. Thus our results show that at least 
$47.94\%$ of all numerical sets in $\Scal(g)$ are irreducible. Another nice property is that 
the type function is multiplicative when restricted to 
$ \bigcup \setpres{{\Gcal}(g)}{$g \ge 1$}$ by \cite[Proposition~5.3]{AM} 
(that is, the type of a
direct sum is the product of the types of its factors, if the factors have no small atoms). Multiplicativity of type was a central theme in \cite{AM}. 
We also mention that when a numerical set $S$ is in $\Gcal(g)$ its type can be computed 
via the formula
  $$ type(S) = \frac{\card{S \cap [0,g)} \, \card{S^* \cap [0,g)} }{\card{S \cap S^* \cap [0,g)}^2} ,$$ 
which is readily derived from, but much simpler than, the general formula for the type of 
an arbitrary numerical set given in \cite{AM}.

In the last two sections of the paper we explore parallel ideas for counting the
number of symmetric numerical sets in $\Gcal(g)$. This study is suggested by Backelin's 
examination
of the number of symmetric numerical monoids in $\Scal(g)$ in \cite{B}. We show that the 
ratio $\gsymm{g}$ of the number of symmetric numerical sets in $\Gcal(g)$ by the total 
number of symmetric numerical sets with Frobenius number $g$ has a limit $\gsymm{\infty}$
which is approximately equal to $.23644$. We also obtain information about the singularities
of the generating function for the sequence $\set{\gsymm{g}}$.
In many ways the analysis of $\set{\gsymm{g}}$ turns out to be more elementary than that
of $\set{\gamma_g}$. For example, in the symmetric setting we obtain two recursively related 
sequences of
integers $\set{\Asigma_k}$ and $\set{\Asp{k}}$,
and the odd terms in the second of these sequences coincides with a well-known sequence 
consisting of the numbers of additive $2$-bases for $k$.

\section*{Numerical sets with no small atoms}

Let $S$ be a numerical set with Frobenius number $g$.
A {\it small atom} for $S$ is a (nonzero) atom for $S$
which is less than $g$. 
A small atom for $S$ is said to be {\it large} if it is greater
than $g/2$. The first result says that every numerical set
which has a small atom will have a large small atom.

\begin{lem}\label{lem:LargeSmallAtom}
Let $S$ be a numerical set in $\Scal(g)$. If $S$ has a small 
atom then $S$ has a small  atom larger than $g/2$.
\end{lem}

\begin{proof} If $g$ is even then $g/2$ is not an atom of $S$ since 
$g/2+g/2 \notin S$. 
Suppose $S$ has an atom less than $g/2$
and let $k$ be the largest such atom.
Then $2k$ is a small atom of $S$, and $2k$ is greater
than $g/2$ by the choice of $k$.       
\end{proof}

The set $\Scal(g)$ is partitioned into two subsets
  $$\Gcal(g)  \ = \ \Gcal(\Nbb_g) \ = \ \setpres{S \in \Scal(g)}{$S$ has no small atoms}$$
and 
  $$\Bcal(g) \ = \ \setpres{S \in \Scal(g)}{$S$ has at least one small atom} \ .$$
For each $g>0$, $\Nbb_g \in \Gcal(g)$ and $\Gcal(g)$ is nonempty.
On the other hand, $\Bcal(g)$ contains all of the numerical monoids in $\Scal(g)$ other than
$\Nbb_g$, and $\Bcal(g)$ is nonempty when $g > 2$.
We are interested in the two ratios
  $${\beta}_g = \frac{\card{\Bcal(g)}}{\card{\Scal(g)}} = \frac{\card{\Bcal(g)}}{2^{g-1}}$$ 
and
   $${\gamma}_g = \frac{\card{\Gcal(g)}}{\card{\Scal(g)}}= \frac{\card{\Gcal(g)}}{2^{g-1}}.$$
Observe that $0 \leq {\beta}_g,{\gamma}_g \leq 1$ and that ${\beta}_g + {\gamma}_g=1$.  

For each $S \in \Scal(2n-1)$ and $\epsilon \in \Zbb_2=\set{0,1}$ we define
   \[ S_\epsilon^\prime = \big(S \cap [0,n-1]\big) \cup  \set{\epsilon \, n} \cup 
      \big(1+S \cap [n,\infty)\big) .\]

\begin{figure}[hbt]
\bigskip
    \centering
   \includegraphics[width=4.0in]{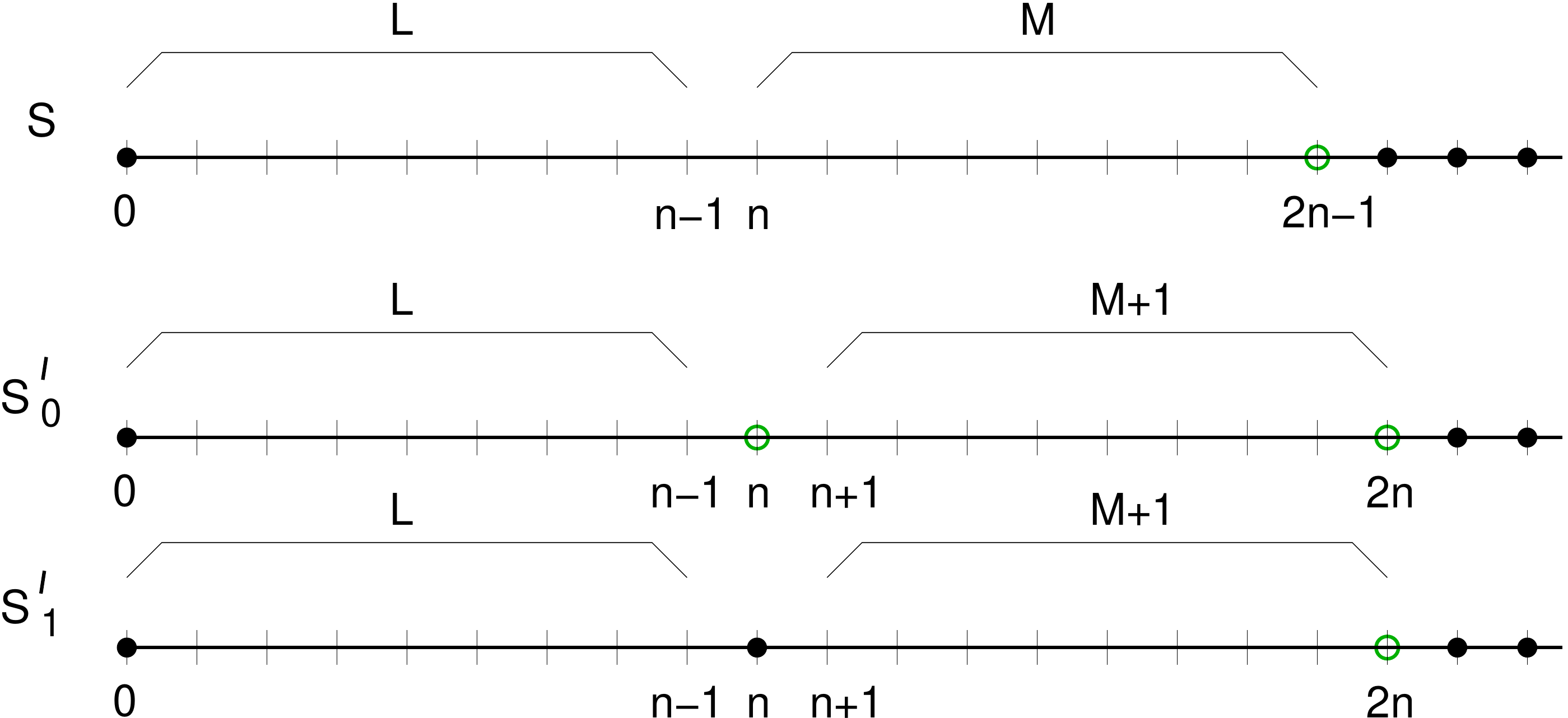}
   \caption{\textsl{The numerical sets $S_0^\prime$ and $S_1^\prime$}}\label{fig:Sepsilon}
\medskip
\end{figure}

\begin{lem}\label{lem:EvenOdd} 
The correspondence $(S,\epsilon) \mapsto S_\epsilon^\prime$ 
is a  bijection from $\Scal(2n-1) \times \Zbb_2$ to $\Scal(2n)$ which
carries $\Gcal(2n-1)\times \Zbb_2$ onto $\Gcal(2n)$. Furthermore, $\gamma_{2n}= 
\gamma_{2n-1}$ and $\beta_{2n}= \beta_{2n-1}$.
\end{lem}

\begin{proof}
The correspondence $(S,\epsilon) \mapsto S_\epsilon^\prime$ is injective by definition,
and it is also surjective: if $S^\prime \in \Scal(2n)$ then $S^\prime=S_\epsilon^\prime$
where $S$ is the union of $S^\prime \cap [0,n-1]$ and $-1+\left( S \cap [n+1,\infty) 
\right)$, and $\epsilon$ equals $0$ if $n \notin S^\prime$ and $1$ if $n \in S^\prime$.

It is not difficult to see that an integer $x$ is a large small atom for $S$ if
and only if $1+x$ is a large small atom for $S_\epsilon^\prime$. By 
lemma~\ref{lem:LargeSmallAtom} this implies that a numerical set $S \in \Scal(2n-1)$
is in $\Gcal(2n-1)$ if and only if $S_\epsilon^\prime$ is in $\Gcal(2n)$. 
To complete the proof, we note that 
$ \gamma_{2n} = \card{\Gcal(2n)}/2^{2n-1} = 2\card{\Gcal(2n-1)}/2^{2n-1} = \gamma_{2n-1}$.
\end{proof}

For integers $g$ and $k$ with $g > k > 0$, let    
	$$\Bcal(g,k) \ = \ \setpres{S \in \Scal(g)}{$g-k$ is the largest small atom of $S$} \ .$$
Note that $\Bcal(g,k)$ is a subset of $\Bcal(g)$ and   
that $\Bcal(g,k)$ is empty whenever $k \ge g/2$ by lemma~\ref{lem:LargeSmallAtom}.
In order to describe $\Bcal(g,k)$ we are led to the next definition.
An ordered pair $(L,M)$ of subsets of $(0,k)$ is {\it admissible} if it satisfies two conditions:
   \begin{itemize}
       \item[(ad1)]  $L \subset M$, and
       \item[(ad2)]  for every $x \in M$ there exists $y \in L$ with
            $x+y \leq k$ and $x+y \notin M$.
   \end{itemize}
Let $\Acal(k)$ be the set of all admissible pairs of subsets of $(0,k)$,
and let $A_k=\card{\Acal(k)}$ denote the cardinality of this set. 
The power set $\Pcal(k,g-k)$ of the set $(k,g-k)$ consists of all subsets of $(k,g-k)$ and has cardinality
$2^{g-2k-1}$.

\begin{thm}\label{thm:CountBcal} 
For integers $g$ and $k$ with $g > 2k > 0$
the set $\Bcal(g,k)$ is in one-to-one correspondence
with ${\Acal(k)} \times \Pcal(k,g-k)$.  In particular, the cardinality of $\Bcal(g,k)$
equals $A_k \, 2^{g-2k-1}$.
\end{thm}

\begin{figure}[htb]
\bigskip
    \centering
   \includegraphics[width=4.0in]{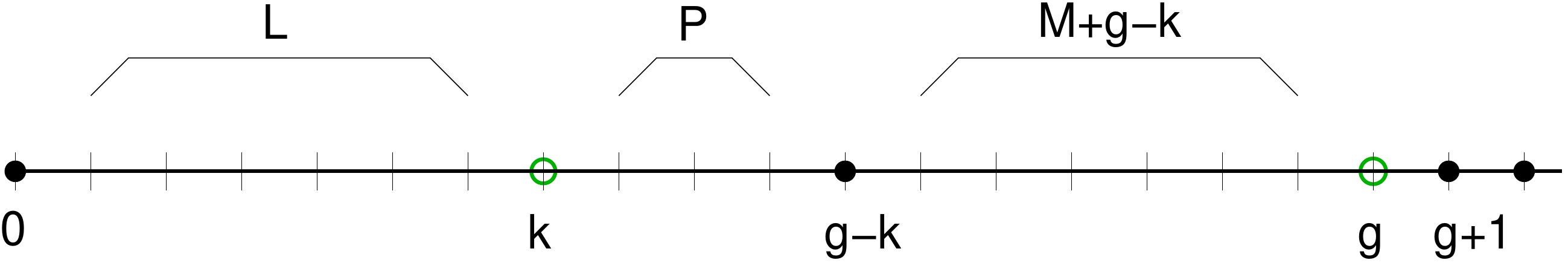}
   \caption{\textsl{The numerical set $S(L,M,P)$}}\label{fig:SLMP}
\medskip
\end{figure}

\begin{proof}
Suppose that $(L,M) \in \Acal(k)$ and  $P \in \Pcal(k,g-k)$. Then
  \begin{equation}\label{eqt:S(L,M,P)}
    S(L,M,P)= \Nbb_g \cup L \cup P \cup \set{g-k} \cup \big(g-k+M\big) 
  \end{equation}    
is a numerical set with Frobenius number $g$. (See figure~\ref{fig:SLMP}.)
Since $g-k+L \subseteq g-k+M$ by (ad1) and each nonzero element of $S(L,M,P)-L$ is larger than $k$, 
$g-k$ is a small atom for $S(L,M,P)$.
Suppose that $x \in (g-k,g) \cap S(L,M,P)$. Then $x-g+k \in M$ and by (ad2) there is an
integer $y \in L$ such that $y+x <g$ and $y+x \notin g-x+M$. This shows that $y+x \notin S(L,M,P)$
and that $x$ is not an atom for $S(L,M,P)$. Thus $g-k$ is the largest small atom for 
$S(L,M,P)$ and $(L,M,P) \mapsto S(L,M,P)$ describes a function $\theta$ from ${\Acal(k)} \times \Pcal(k,g-k)$
into $\Bcal(g,k)$.

Now assume that $S\in \Bcal(g,k)$, and define $L_S \subseteq (0,k)$, $M_S \subseteq (0,k)$ and 
$P_S \subseteq (g-k,g)$ by 
  $$L_S=S \cap (0,k), \ M_S = k-g+ \big( S \cap (g-k,g) \big) , \ 
      P_S=S \cap (k,g+k) \ .$$  
Since $g-k$ is a small atom of $S$ then $\ell +g-k \in S \cap (g-k,g)$ for each 
$\ell \in L_S$, which implies that $\ell \in  M_S$ and that $L_S \subseteq M_S$.
Suppose that $x \in M_S$. Then $g-k+x$ is an element of $S$ but not an atom
of $S$ (since $g-k$ is the largest small atom of  $S$) and so  
there exists $y \in S$ with $g-k+x+y \le g$ and 
$g-k+x+y \notin S$. It follows that $x+y \le k$,
$y \in L_S$ and $x+y \notin M_S$.
Thus the pair $(L_S,M_S)$ satisfies (ad1) and (ad2) and $(L_S,M_S) \in \Acal(k)$.
Let $\Phi$ be the function from $\Bcal(g,k)$ to  
${\Acal(k)} \times \Pcal(k,g-k)$ given by $S \mapsto (L_S,M_S,P_S)$.
The proof is completed by observing that $\theta$ and $\Phi$ are inverses of each other. 
\end{proof}

By lemma~\ref{lem:LargeSmallAtom} the set $\Bcal(g)$ can be expressed as the
disjoint union of the sets $\Bcal(g,k)$ where $k$ ranges from $1$ to $\floor{(g-1)/2}$.
Thus, the cardinality of $\Bcal(g)$ is the
sum of the cardinalities of $\Bcal(g,k)$ for $1 \le k \le \floor{(g-1)/2}$.
	
\begin{cor}\label{cor:betaSum}
For each positive integer $g$, ${\displaystyle {\beta}_g = \sum_{k=1}^{\floor{(g-1)/2}} A_k \, 4^{-k}}$. 
\end{cor}

\begin{proof} \
Because of theorem~\ref{thm:CountBcal} and the comment above, we have
 $$ \beta_g = \frac{\card{\Bcal(g)}}{2^{g-1}} = \sum_{k=1}^{\floor{(g-1)/2}} \frac{\card{\Bcal(g,k)}}{{2^{g-1}}}  
    = \sum_{k=1}^{\floor{(g-1)/2}} \frac{A_k \, 2^{g-2k-1}}{2^{g-1}} = \sum_{k=1}^{\floor{(g-1)/2}} A_k \, 4^{-k} .$$   	
\end{proof}

By corollary~\ref{cor:betaSum} the sequence $\set{\beta_g}$ is increasing, and
being bounded above by $1$, it must have a limit
  $${\beta}_{\infty} = \lim_{g \to \infty} {\beta}_g = \sum_{k=1}^{\infty}
	   A_k 4^{-k} \ .$$
As a consequence the sequence $\set{\gamma_g}=\set{1 - {\beta}_g}$ is decreasing 
with limit
  $${\gamma}_{\infty} = \lim_{g \to \infty} {\gamma}_g = 1 - {\beta}_{\infty} \ .$$
By the next lemma, it is also possible to express $\gamma_\infty$ as the sum 
of a positive series ${\gamma}_{\infty}=\sum_{k=1}^{\infty} (3^{k-1}-A_k) 4^{-k}$.

\smallskip

\begin{lem}\label{lem:A_kBounds}
For each integer $k>0$, $2^{\floor{(k-1)/2}} \leq A_k \leq 3^{k-1}$. 
Moreover $\gamma_{2k-1} - \gamma_\infty$ is positive and
$\gamma_{2k-1} - \gamma_\infty =  \beta_\infty - \beta_{2k-1}  
        \le \left({3}/{4}\right)^{k-1}$. 
\end{lem}

\begin{proof}
Let $L$ be an arbitrary nonempty subset of $(0,\floor{(k+1)/2}) \subset (0,k)$
with maximal element  $\ell$.  For any element
$x \in L$, $\ell+x \leq 2\ell \leq k$ and $\ell+x \notin L$.  This shows that $(L,L)$ 
is an admissible pair of subsets of $(0,k)$.
Since there are $2^{\floor{(k-1)/2}}$ distinct subsets of  $(0,\floor{(k+1)/2})$,
this verifies the inequality $2^{\floor{(k-1)/2}} \leq A_k$. 

Suppose $(L,M)$ is an element of $\Acal(k)$. Then for each $x \in (0,k)$ we have three distinct
possibilities: (1) $x \notin M$, (2) $x \in M$ and $x \notin L$, or (3) $x \in L$.
Therefore there are $3^{k-1}$ pairs of subsets $(L,M)$ in $(0,k)$ which satisfy
$(ad1)$, and it follows that $A_k \le 3^{k-1}$. Now by definition and corollary~\ref{cor:betaSum} 
 $$\gamma_{2k-1} - \gamma_\infty =  \beta_\infty - \beta_{2k-1}
   = \sum_{i=k}^\infty \, A_i \, 4^{-i} \le  \sum_{i=k}^\infty \, 3^{i-1} \, 4^{-i} = (3/4)^{k-1}
   \ .$$  
\end{proof}

Notice that $(\emptyset,\emptyset)$ is the only ordered pair of subsets of $(0,1)=\emptyset$, and as
it is admissible, this shows that $A_1=1$. Among ordered pairs of subsets of $(0,2)$, condition
(ad1) fails for $\left(\set{1},\emptyset\right)$ and condition (ad2) fails for $(\emptyset,\set{1})$ 
while the two remaining ordered pairs $(\emptyset,\emptyset)$ and $(\set{1},\set{1})$ are in 
$\Acal(2)$, and so $A_2=2$.  With lemma~\ref{lem:A_kBounds} and these values of $A_1$ and $A_2$,
  $$ \beta_\infty \le \beta_5 + (3/4)^2 = \left( \frac{1}{4} + \frac{2}{16}\right) + \frac{9}{16} 
           = \frac{15}{16} ,$$
which shows that both $\gamma_\infty$ and $\beta_\infty$ are strictly between $0$ and $1$.
Using this approach with the more extensive data compiled in table~\ref{tab:beta}, we see that 
$\beta_{33} = .510538\ldots$ approximates $\beta_\infty$ to within $(3/4)^{16} = .0100226\ldots$.
Taking midpoints gives the approximation $\beta_\infty \approx .515549$ accurate to within
$.005011$, and subtracting from $1$ gives $\gamma_\infty  \approx .484451$ 
with the same degree of accuracy. This approximation can be rephrased as saying that 
$\card{\Gcal(g)} \approx .484451 \times 2^{g-1}$ 
for large values of $g$.

\begin{table}[ht]  
\bigskip
\begin{tabular}{|c||c|c|c|c|c|}
\hline
 $n$ & $A_n^\prime$ & $A_n$ & $\beta_{2n+1}$ & $\beta_{2n+1} + (3/4)^n$ & $A_{n-1}/A_n$\\
\hline
   1  & 1     & 1       & \, .250000 \, & 1.000000 & -\\
   2  & 2     & 2       & .375000 & .937500  & .5000\\
   3  & 3     & 3       & .421875 & .843750  & .6667\\
   4  & 6     & 8       & .453125 & .769531 & .3750\\
   5  & 10    & 18      & .470703 & .708008 & .4444\\
   6  & 20    & 50      & .482910 & .660889 & .3600\\
   7  & 37    & 135     & .491150 & .624634 & .3704\\
   8  & 74    & 385     & .497025 & .597137 & .3506\\
   9  & 140   & 1065    & .501087 & .576172 & .3615\\
   10 & 280   & 3053    & .503999 & .560312 & .3488\\
   11 & 542   & 8701    & .506073 & .548308 & .3509\\
   12 & 1084  & 25579   & .507598 & .539274 & .3402\\
   13 & 2118  & 73693   & .508696 & .532453 & .3471\\
   14 & 4236  & 217718  & .509507 & .527325 & .3385\\
   15 & 8337  & 635220  & .510090 & .523462 & .3427\\
   16 & \, 16674 \, & \, 1888802 \, & .510538 & .520561 & .3363\\
\hline
  \end{tabular}
\vskip6pt
\caption{Bounds for ${\beta}_{\infty}$.}\label{tab:beta}
\end{table}

If $(L,M)$ is an admissible ordered pair of subsets of $(0,k)$ and
$L^\prime$ and $M^\prime$ are subsets satisfying $L \subseteq 
L^\prime \subseteq M^\prime \subseteq M$ then $(L^\prime,M^\prime)$
is also admissible. In particular, both $(L,L)$ and $(M,M)$ are
elements of $\Acal(k)$ whenever $(L,M) \in \Acal(k)$. The computer
routine that was used to generate the data in table~\ref{tab:beta}
starts by first determining the collection of 
subsets $L \subseteq (0,k)$ for which $(L,L)$ is admissible. (The cardinality
of this collection is denoted by $A_k^\prime$ in the table.
These numbers are important in their own right as we shall explain in 
the next section.)
The routine then isolates nested pairs of sets in this collection and
tests only these pairs for condition (ad2). Even with this, the algorithm
has exponential complexity and slows down quite rapidly.

\medskip

From the results of this section we may draw some further conclusions 
which directly address the Anti-Atom Problem for an arbitrary numerical monoid $M$.

\begin{thm}\label{thm:AtomSet} Let $M \ne \Nbb_g$ be a numerical monoid with 
Frobenius number $g$ and let $g-k$ be the largest element of $M \cap (0,g)$. Then
$\card{\Gcal(M)}\le A_k 4^{-k}\, 2^{g-1}\le \frac{1}{3}\left( 3/4 \right)^k \times 2^{g-1}$.
\end{thm}

\begin{proof}
If $M \ne \Nbb_g$ is a numerical monoid in $\Scal(g)$ and $g-k$ is the largest 
element in $M \cap (0,g)$ then $g-k$ is the largest small atom of every numerical set $S$ with $A(S)=M$. 
Thus $\Gcal(M) \subseteq \Bcal(g,k)$ and
$\card{\Gcal(M)} \le \card{\Bcal(g,k)} = A_k 2^{g-2k-1}$. The last inequality follows
from lemma~\ref{lem:A_kBounds}. 
\end{proof}

The value of $k$ in theorem~\ref{thm:AtomSet} satisfies $0 < k < g/2$. Since 
$\frac{1}{3}\left( 3/4 \right)^k \times 2^{g-1}\le .25 \times 2^{g-1}$ is less than 
$.484451 \times 2^{g-1}$ for all $k$, we see that among all monoids in $\Scal(g)$ the 
one with largest anti-atom set is always $\Nbb_g$ (which is not too surprising 
since $\Gcal(\Nbb_g)$ contains more than $48\%$ of the elements of $\Scal(g)$).

As $k$ increases from $1$ to $\floor{(g-1)/2}$ the cardinality of $\Bcal(g,k)$
decreases but the number of monoids in $\Bcal(g,k)$ decreases as well. For example,
$\Bcal(g,1)$ contains all of the symmetric and pseudosymmetric monoids in $\Scal(g)$,
while, at the other extreme, $\Bcal(g,\floor{(g-1)/2})$ contains only one monoid $\Dbb_g$
which is defined by
   \begin{equation}\label{eqt:Dbb}
     \Dbb_g=\Nbb_g \cup {\set{\floor{(g+2)/2}}} .
   \end{equation}

\smallskip

\begin{cor}
For each nonnegative integer $n$, we have $\card{\Gcal(\Dbb_{2n+1})}= A_n$ and $\card{\Gcal(\Dbb_{2n+2})}= 2A_{n}$.
 \end{cor} 

\begin{proof}
It is not difficult to show  that $\Dbb_g$ is the only monoid in $\Bcal(g,\floor{(g-1)/2})$.
(If $M$ is a monoid in $\Bcal(g,\floor{(g-1)/2})$ then $M \cap (\floor{(g+2)/2}, g) = \emptyset$.)
Thus $\Dbb_g$ is the atom monoid of every numerical set in $\Bcal(g,\floor{(g-1)/2})$,
and this implies that $\card{\Gcal(\Dbb_g)} =
\card{\Bcal(g,\floor{(g-1)/2})} = A_{\floor{(g-1)/2}} 2^{g-2\floor{(g-1)/2}-1}$,
from which the corollary follows.
\end{proof}

When $k$ is less than $\floor{(g-1)/2}$ the set $\Bcal(g,k)$ will always contain at least
two distinct numerical monoids (for example, $\Nbb_g \cup \set{g-k}$ and $\Nbb_g \cup \set{g-k-1,g-k}$).
Thus $\Dbb_g$ is the only monoid in $\Scal(g)$ for which the first inequality of 
theorem~\ref{thm:AtomSet} is sharp.

\section*{The generating function for $\set{\gamma_k}$}

For each integer $k>0$ let $\Acal(k)^\prime$ denote the collection of
all subsets $L \subseteq (0,k)$ for which $(L,L)$ is admissible. 
Thus $\Acal(k)^\prime$ consists of those subsets $L$ which satisfy 
the condition that
for each $x \in L$ there is $y \in L$ such that $x+y \le k$ and $x+y \notin L$.
The cardinality of $\Acal(k)^\prime$ will be denoted by $A_k^\prime$.

\begin{thm}\label{thm:GoodSet}
There is a one-to-one correspondence between $\Gcal(g)$ and $\Acal(g)^\prime$.  
In particular, $\card{\Gcal(g)}=A_g^\prime$ and $\gamma_g = A_g^\prime \, / 2^{g-1}$.
\end{thm}

\begin{proof} 
For each $L \in \Acal(g)^\prime$ consider the numerical set $\Nbb_g \cup L \in \Scal(g)$.
If $x \in L$ then there is an integer $y \in L$ such that 
$0 < x+y \le g$ and $x+y \notin L$, which implies that $x+y \notin \Nbb_g\cup L$. 
This shows that $x$ is not an atom of $\Nbb_g \cup L$, and that 
$\Nbb_g \cup L$ has no small atoms. 
Thus the correspondence $L \mapsto \Nbb_g \cup L$ is a function from $\Acal(g)^\prime$ to
$\Gcal(g)$, and clearly this function is injective.
Now suppose $S \in \Gcal(g)$ and let  $x \in S\cap (0,g)$. 
Since $S$ has no small atoms, there is an integer $y \in S$ such that $x+y \notin S$.
Thus $x+y \le g$ (since the Frobenius number of $S$ is $g$),
$y \in S \cap (0,g)$ and $x+y \notin S \cap (0,g)$. 
It follows that $S \cap (0,g)$ is an element of $\Acal(g)^\prime$, and the function 
$L \mapsto \Nbb_g \cup L$  is surjective. 
\end{proof}

\begin{thm}\label{thm:recursions}
For each $k \geq 1$, $A_{2k}^{\prime}=2 \ A_{2k-1}^{\prime}$ and
$A_{2k+1}^{\prime}=2 \ A_{2k}^{\prime} - A_k$.
\end{thm}

\begin{proof} 
The first equation follows immediately from lemma~\ref{lem:EvenOdd} and
theorem~\ref{thm:GoodSet}. For the second equation we have
 $$A_{2k+1}^\prime = \card{\Gcal(2k+1)} = \card{\Scal(2k+1)} - \card{\Bcal(2k+1)}
    = 2^{2k} - \sum_{\ell =1}^k \card{\Bcal(2k+1,\ell)} ,$$  
and by theorem~\ref{thm:CountBcal}
 $$ 2^{2k} - \sum_{\ell =1}^k \card{\Bcal(2k+1,\ell} = 4^k - \sum_{\ell =1}^k A_\ell \, 4^{k-\ell} \ .$$
A similar computation shows that
 $$A_{2k}^\prime = \frac{1}{2}\Big(4^k - \sum_{\ell =1}^{k-1} A_\ell \, 4^{k-\ell}\Big) \ .$$
Combining these gives $A_{2k+1}^\prime -  2A_{2k}^\prime = -A_k$.  
\end{proof}

\begin{figure}[hbt]
\bigskip
    \centering
   \includegraphics[width=4.9in]{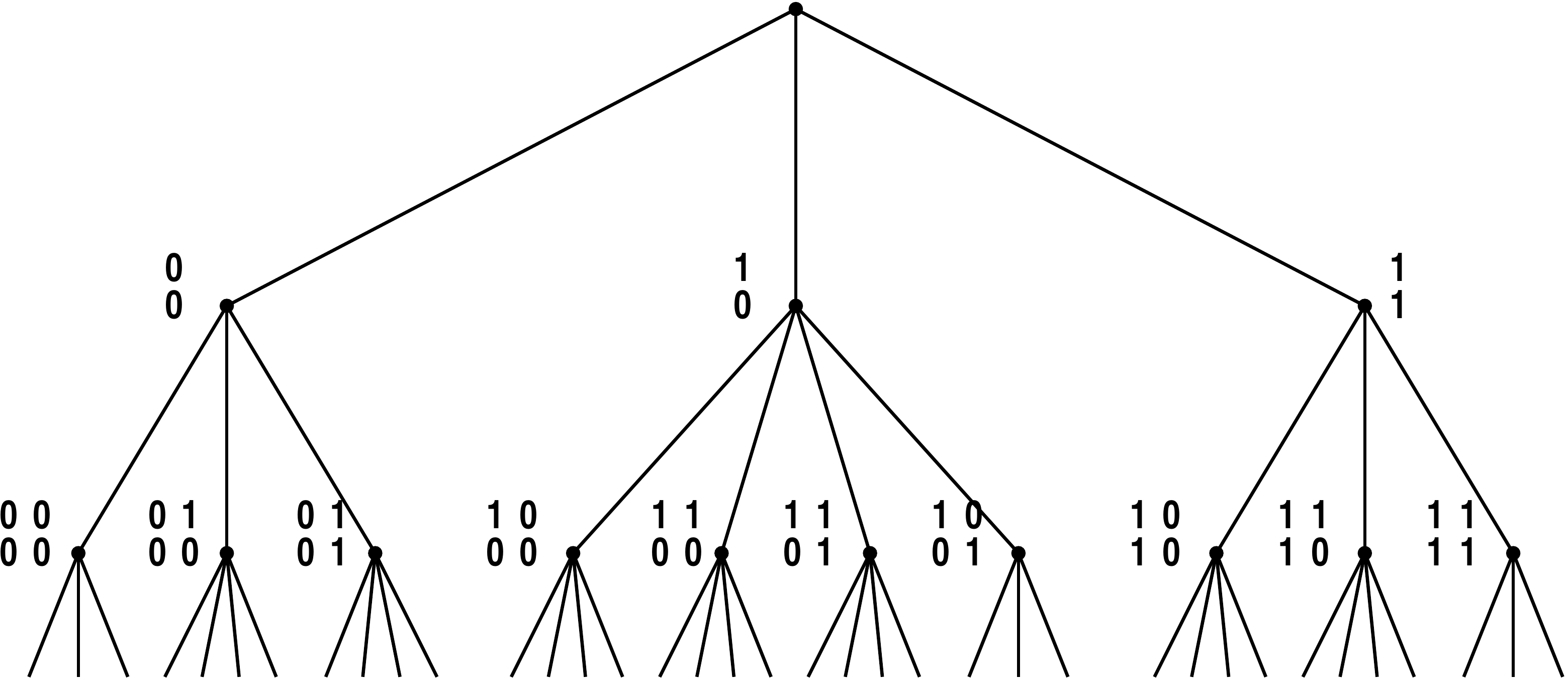}
   \caption{\textsl{Rooted tree for \ $\bigcup \setpres{\Gcal(2k+1)}{$k\in                            \Nbb$}$.}}\label{fig:RootedTree}
\medskip
\end{figure}

The set $\Gcal(2k+1)$ can be constructed from $\Gcal(2k-1)$ by a process in which
each element of $\Gcal(2k-1)$ will spawn either three or four elements of $\Gcal(2k+1)$
as follows.  If $S \in \Gcal(2k-1)$ and $Q$ is one of the four subsets of $\set{k,k+1}$
then let $S(Q) \in \Scal(2k+1)$ be given by 
   \begin{equation}\label{eqt:spawn} 
      S(Q) = \Nbb_{2k+1} \cup \Big( S\cap [1,k-1] \Big)  \cup   Q  \cup  \Big( 2+\big(S\cap [k,2k-2]\big) \Big) \ .
   \end{equation}
Since $S$ has no small atoms, $S(Q)$ will not have any small atoms larger than $k+1$.
\begin{figure}[hbt]
\bigskip
    \centering
   \includegraphics[width=4.0in]{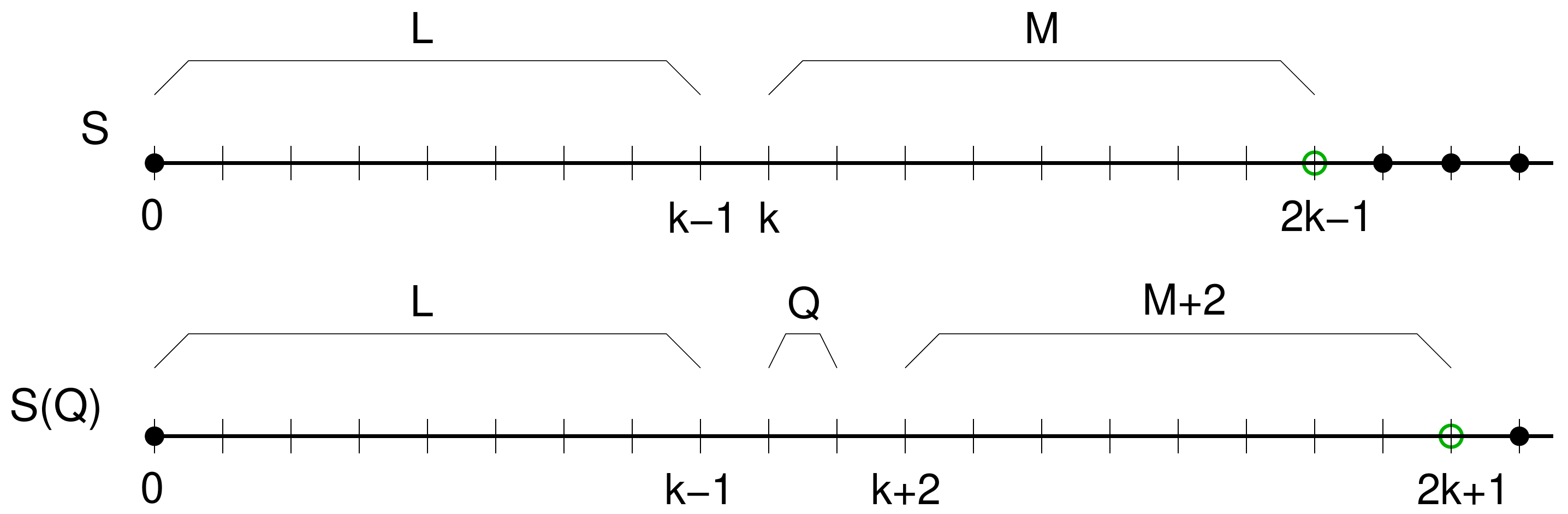}
   \caption{\textsl{The numerical set $S(Q)$}}\label{fig:SQ}
\medskip
\end{figure}
Furthermore, if $Q$ is one of $\emptyset$, $\set{k}$ or $\set{k,k+1}$ then $k+1$ is not an atom of $S(Q)$, 
and $S(Q)$ is an element of $\Gcal(2k+1)$.
But if $Q$ equals the singleton set $\set{k+1}$ then sometimes $k+1$ will be an atom for $S(Q)$, 
in which case $S(Q)$ is not an element of $\Gcal(2k+1)$. 
From this we see that $\card{\Gcal(2k+1)} = A^\prime_{2k+1}$ is four times $\card{\Gcal(2k-1)}=A^\prime_{2k-1}$
minus the number of elements of  $\Gcal(2k-1)$ which spawn only three elements of $\Gcal(2k+1)$, and,
since $A^\prime_{2k+1} = 4A^\prime_{2k-1} - A_k$ by  theorem~\ref{thm:recursions}, 
the number of elements of  $\Gcal(2k-1)$ which spawn only three elements of $\Gcal(2k+1)$ equals 
$A_k$.\footnote{This can also be verified by a combinatorial argument.
If $(L,M) \in \Acal(k)$ then $(L^\prime,L^\prime)$
is an element of $\Acal(2k-1)^\prime \cong \Gcal(2k-1)$ which spawns only three elements of
$\Gcal(2k+1)$, where $L^\prime = L \cup \left( M+k-1 \right)$.}      
As a result of these comments we can view the union of all the sets $\Gcal(2k+1)$ 
as the vertices of a downward opening rooted tree in which each vertex is directly
above the $3$ or $4$ vertices that it spawns,  as pictured in figure~\ref{fig:RootedTree}. In the illustration 
the vertex labeled by a $2 \times k$ matrix $\alpha = \left( \begin{smallmatrix} \alpha_1 & \alpha_2 &\cdots &\alpha_k 
\\ \alpha_{2k} & \alpha_{2k-1} &\cdots &\alpha_{k+1} \end{smallmatrix} \right)$ 
with entries in $\Zbb_2$ corresponds to the numerical set 
  $$ S(\alpha) = \Nbb_{2k+1} \cup \setpres{i}{$\alpha_i=1$} \ $$
in $\Gcal(2k+1)$.  
Although we will not use it here, one can specify conditions on the matrix $\alpha$ which 
guarantee that $S(\alpha)$ is in
$\Gcal(2k+1)$: Call the $2 \times k$ matrix $\alpha$ {\it quadrivalent}  if there is
an integer $\ell$ with $1 \le \ell \le k$ such that the $\ell$th column of $\alpha$ is 
$\left( \begin{smallmatrix} 1\\ 0 \end{smallmatrix} \right)$  or $\left( \begin{smallmatrix} 1\\ 1 \end{smallmatrix} \right)$
and the $(k+1-\ell)$th column is $\left( \begin{smallmatrix} 1\\ 0 \end{smallmatrix} \right)$  or $\left( \begin{smallmatrix} 0\\ 0 \end{smallmatrix} \right)$.\footnote{The reason for this terminology is that if $S(\alpha)$ 
is an element of $\Gcal(2k+1)$, then the matrix $\alpha$ is quadrivalent if and only if $S(\alpha)$
spawns four elements of $\Gcal(2k+3)$.}
Then $S(\alpha) \in \Gcal(2k+1)$ if and only if whenever the $i$th column of
$\alpha$ equals $\left( \begin{smallmatrix} 0\\ 1 \end{smallmatrix} \right)$ then the $2 \times (i-1)$ submatrix
of $\alpha$ to the left of that column is quadrivalent. 
We also note that as one moves down the tree the ratio $A_k/A^\prime_{2k+1}$ of the number of 
vertices at level $2k+1$ which spawn three vertices by the total number of vertices at that level limits to $0$.
Indeed  $A_k$ is bounded above by the number of $2 \times k$ matrices $\alpha$ which are not quadrivalent 
and that number is easily seen to equal $3^k$. Thus 
 $$\frac{A_k}{A^\prime_{2k+1}} = \frac{A_k}{\abs{\Gcal(2k+1)}} = \frac{A_k}{\gamma_{2k+1}4^k} \le 
   \frac{1}{\gamma_{2k+1}}\left(\frac{3}{4}\right)^k$$  
and the latter limits to $0$.

Let $g(z)$ and $f(z)$ be the analytic functions defined by 
$$   g(z) = \sum_{k=1}^\infty \, A_k z^k
 \qquad \text{and} \qquad  
       f(z) = \sum_{k=1}^\infty \, A_k^\prime z^k \ .$$
  
\begin{cor}\label{cor:fgrelation} 
The functions $f(z)$ and $g(z)$ satisfy the relation  
  $$\big(2z-1\big)f(z) = z\left(g(z^2)-1\right) \ .$$
\end{cor}

\begin{proof}
Using theorem~\ref{thm:recursions} we have
\begin{equation}\label{eqt:f(z)}
f(z) =  \sum_{k=1}^\infty \, A_{2k-1}^\prime z^{2k-1} +  \sum_{k=1}^\infty \, A_{2k}^\prime z^{2k}
     = (2z+1) \sum_{k=1}^\infty \, A_{2k-1}^\prime z^{2k-1} \ .
\end{equation} 
Also $A_k = 4A_{2k-1}^\prime - A_{2k+1}^\prime$ by theorem~\ref{thm:recursions}, and  
$$\begin{aligned} g(z^2) &=  
  \sum_{k=1}^\infty \, (4A_{2k-1}^\prime - A_{2k+1}^\prime)z^{2k}\cr 
    &=  4z\sum_{k=1}^\infty \, A_{2k-1}^\prime z^{2k-1} -
    \frac{1}{z}\sum_{k=1}^\infty \, A_{2k-1}^\prime z^{2k-1} +   A_{1}^\prime     \cr
    & = \frac{(2z - 1)(2z+1)}{z} \sum_{k=1}^\infty \, A_{2k-1}^\prime z^{2k-1} + 1
    = \frac{2z-1}{z}f(z) + 1\  ,
    \end{aligned}$$ 
where the last equality follows from (\ref{eqt:f(z)}).  
\end{proof}

\begin{cor}\label{cor:f(z)Radius}  
The analytic function $f(z)$ has a 
singularity at $z=1/2$, and its radius of convergence at the origin equals $1/2$. 
Other than $z=1/2$, the singularities of $f(z)$ coincide with those of $g(z^2)$
and $f(z)(z-1/2)$ is continuous on $\abs{z} \le 1/2$.
\end{cor}

\begin{proof} 
Since $ \sum_{k=1}^\infty \, A_k 4^{-k}$ sums to $\beta_\infty = g(1/4)$, 
$\sum_{k=1}^\infty \, A_k z^k$ converges absolutely and $f(z)$ is
continuous on $\abs{z} \le 1/4$. 
Note that $f(z)$ has no singularities in $\abs{z} < 1/2$ because $g(z)$ has
none in $\abs{z} < 1/4$.  The remainder of the proof follows from corollary~\ref{cor:fgrelation}.
\end{proof}

By definition the integers $A_k^\prime$ and $A_k$ satisfy $0 \le A_k^\prime \le A_k$ for all $k>0$, and by corollary~\ref{cor:f(z)Radius} we know that
the series $\sum_{k=1}^\infty \, A_k^\prime z^k$ diverges when $z=1/2$.
Therefore $\sum_{k=1}^\infty \, A_k z^k$ must also diverge for $z=1/2$ by the 
comparison test, and this 
shows that the radius of convergence for $g(z)$ at the origin is between $1/4$ and $1/2$. To find 
the precise value, the ratio test would lead one to examine the sequence $A_{n-1}/A_n$. Empirical evidence from the last column of table~\ref{tab:beta} perhaps suggests that this sequence
has a limit infimum larger than $1/4$, but we have not been able to 
ascertain this. So the determination of the radius of convergence of $g(z)$ 
at the origin remains open.
Notice that if this radius of convergence is larger than $1/4$ then
$g(z^2)$ is analytic in a neighborhood of $z=1/2$ and
  $$ \lim_{z\to 1/2}\,(z-1/2)f(z) =\lim_{z\to 1/2}\,{z}\left(g(z^2) -1\right)/2 
     = (\beta_\infty - 1)/4= -\gamma_\infty/4 ,$$
which would imply that $f(z)$ has a simple pole with residue $-\gamma_\infty/4$ at $z=1/2$.  

\medskip

Now the generating function for the sequence $\set{\gamma_k}$ is 
   $$h(z) =  \sum_{k=1}^\infty \, \gamma_k \, z^k$$
which satisfies
  $$ h(z)  = \sum_{k=1}^\infty \, \frac{A_k^\prime}{2^{k-1}}\, z^k =
    2f\left(z/2\right) = \left( \frac{z}{z-1}\right) \left(g(z^2/4)-1\right) . $$
By corollary~\ref{cor:f(z)Radius} the radius of convergence of this series equals $1$. If
the radius of convergence at the origin for $g(z)$ is  larger than $1/4$ then  
$g(z^2/4)-1$ is analytic in a disk with radius larger than $1$ centered at the origin, and
$h(z)$ has a simple pole at $z=1$ whose residue is $-\gamma_\infty$
(and this is the only pole on $\abs{z}=1$).

\bigskip

\section*{Symmetric numerical sets with no small atoms}

A numerical set $S$ with Frobenius number $g$ is  {\it negative semisymmetric} 
provided that $g-x \notin S$ whenever $x \in S$.
For a positive integer $g$, define
	$$\Ssymm(g) \ = \ \setpres{S \in \Scal(g)}{$S$ is maximal negative semisymmetric 
	   in $\Scal(g)$}$$
where maximality is measured with respect to subset inclusion. 
Then $\Ssymm(g)$ consists of the symmetric
numerical sets in $\Scal(g)$ when $g$ is odd, and the pseudosymmetric numerical sets 
in $\Scal(g)$ when $g$ is even.
Each element of $\Ssymm(g)$ is the union of $\Nbb_g$ 
with a subset of $(0,g)-\set{g/2}$ that is carried onto its complement by $x \mapsto g-x$, and
hence the cardinality of $\Ssymm(g)$ equals $2^{\floor{(g-1)/2}}$.

The set $\Ssymm(g)$ is partitioned into two subsets
  $$\Gsymm(g) \ = \ \setpres{S \in \Ssymm(g)}{$S$ has no small atoms} \ = \ \Gcal(g) \cap \Ssymm(g)$$
and 
  $$\Bsymm(g) \ = \ \setpres{S \in \Ssymm(g)}{$S$ has at least one small atom}  \ = \ \Bcal(g) \cap \Ssymm(g)
  \ .$$
We define  
  $$\bsymm{g} = \frac{\card{\Bsymm(g)}}{\card{\Ssymm(g)}}=\frac{\card{\Bsymm(g)}}{2^{\floor{(g-1)/2}}}$$ 
and
  $$\gsymm{g} = \frac{\card{\Gsymm(g)}}{\card{\Ssymm(g)}}=\frac{\card{\Gsymm(g)}}{2^{\floor{(g-1)/2}}}.$$
The next lemma describes a direct connection between symmetric and pseudosymmetric numerical
sets.

\begin{lem}\label{lem:EvenOddSymm} The correspondence $S \mapsto S_0^\prime$ where
  $$S_0^\prime = \Big(S \cap [0,n-1]\Big) \cup \Big(1+S \cap [n,\infty)\Big)$$ 
defines a bijection from $\Ssymm(2n-1)$ to $\Ssymm(2n)$ 
which carries $\Gsymm(2n-1)$ onto $\Gsymm(2n)$. 
Therefore $\gsymm{2n}= \gsymm{2n-1}$ and $\bsymm{2n}= \bsymm{2n-1}$.  
\end{lem}

\begin{proof} Notice that $S_0^\prime$ is the numerical set $S_\epsilon^\prime$ where 
$\epsilon=0$ as defined in lemma~\ref{lem:EvenOdd}. The proof follows upon observing that
$S$ is symmetric if and only if $S_0^\prime$ is pseudosymmetric.
\end{proof}

For integers $g$ and $k$ with $g > 2k > 0$, let 
	$$\Bsymm(g,k) =  \setpres{S \in \Bsymm(g)}{$g-k$ is the largest small atom of $S$} $$
(i.e.~$\Bsymm(g,k) = \Bcal(g,k) \cap \Ssymm(g)$).
Then $\Bsymm(g)$ is the disjoint union of the sets $\Bsymm(g,k)$ as $k$ ranges between $1$ and 
$\floor{(g-1)/2}$ by lemma~\ref{lem:LargeSmallAtom}. 
In order to describe $\Bsymm(g,k)$ we are led to the next definitions. 

Let $M$ be a subset of $(0,k)$. Define
  \[M_+ = \setpres{m \in M}{$k-m \in M$}  \ ,\]
  \[M_- = \setpres{x \in (0,k)}{$x \notin M$ and $k-x \notin M$}  \ ,\]
and 
  \[M^* = \setpres{x \in (0,k)}{$k-x \notin M$} \ .\footnote{This last definition is closely
  related to the definition of the dual $S^*$ of a numerical set $S$. If $S =\Nbb_k \cup M \in \Scal(k)$ 
  then $S^* = \Nbb_k \cup M^*$.} \]
With these definitions, observe that $M_* = \left( M - M_+\right) \cup M_-$.

A subset $M \subseteq (0,k)$ is called {\it $\sigma$-admissible} if it
satisfies the two conditions:
  \begin{itemize}
       \item[($\sigma$-ad1)] \ $M_- = \emptyset$, and 
       \item[($\sigma$-ad2)] \ for each $x \in M$ there is $y \in M^*$ with
            $x+y < k$ and $x+y \notin M$.
   \end{itemize}
   
Let $\Asymm(k)$ be the set of all $\sigma$-admissible subsets of $(0,k)$ 
and let $\Asigma_k$ denote the cardinality  of this set.  
Also, for integers $k$ and $g$ with $g > 2k > 0$ let $\Psymm(k,g-k)$ be 
the collection of all subsets of $(k,g-k)-\set{g/2}$ that are carried onto 
their complement by the reflection $x \mapsto g-x$. The cardinality of this collection is $2^{\floor{(g-2k-1)/2}}$.
Recall the definition of the numerical set $S(L,M,P)$ as
  $$ S(L,M,P) = \Nbb_g \cup L \cup P \cup \set{g-k} \cup \big(M+g-k\big)$$ 
from equation~(\ref{eqt:S(L,M,P)}) in theorem~\ref{thm:CountBcal}.

\begin{thm}\label{thm:CountBsymm} 
Let $g$ and $k$ be integers with $g > 2k >0$.
The correspondence $(M,P) \mapsto S(M^*,M,P)$ defines a bijection from 
$\Asymm(k) \times \Psymm(k,g-k)$ onto $\Bsymm(g,k)$.
In particular, the cardinality of $\Bsymm(g,k)$
is $\Asigma_k \, 2^{\floor{(g-2k-1)/2}}$.
\end{thm}

\begin{proof} 
Let $(M,P) \in \Asymm(k) \times \Psymm(k,g-k)$. 
An integer $x \in (0,k)$ is an element of $M^*$ if and only if
$k-x \notin M$, which is equivalent to asserting that $g-x \notin
M+g-k$. This together with the fact that $P$ is an element of $\Psymm(k,g-k)$
implies that $S(M^*,M,P) \in \Ssymm(g)$.
Since $M_- = \emptyset$, $M^* = M-M_+ \subseteq M$. 
If $x \in M^*$ then $x \in M$ and $g-k+x \in M+g-k \subset S(M^*,M,P)$, and if 
$x$ is an integer larger than $k$ then $g-k+x > g$ and
$g-k+x \in S(M^*,M,P)$. This shows that
$g-k$ is an atom for $S(M^*,M,P)$. An element of
$S(M^*,M,P)$ in the interval $(g-k,k)$ has the form $x+g-k$ for
some $x \in M$. 
By the definition of $\Asymm(k)$ there is $y \in M^*$ with $y<k-x$
and $x+y \notin M$. 
Thus  $y+(x+g-k)$ is an element of $(g-k,g)$
which is not an element of $M+g-k$, and $x+g-k$ is not an atom of $S(M^*,M,P)$.
It follows that $g-k$ is the largest small atom of $S(M^*,M,P)$ and $S(M^*,M,P) \in \Bsymm(g,k)$.

To complete the proof it only remains to show that each numerical set in  $\Bsymm(g,k)$
equals $S(M^*,M,P)$ for some $(M,P) \in \Asymm(k) \times \Psymm(k,g-k)$. 
Let $T \in \Bsymm(g,k)$, and set $M=\left( T\cap (g-k,g) \right) - g+k \subseteq (0,k)$ and 
$P = T \cap (k,g-k)$.  
If $x \in T\cap (0,k)$ then $x+g-k \in T \cap
(g-k,g)$, since $g-k$ is an atom of $T$, and hence $x \in M$.
Note further that $k-x$ is not in $T$ because if it were then both
$k-x$ and $g-(k-x) = x+g-k$ would be elements of $T$ contradicting
the negative semisymmetry of $T$. This shows that $x \in M^*$ and $T \cap (0,k) 
\subseteq M^*$.
Moreover if $x \in M^*$ then $k-x \notin M$ which means that $g-k+(k-x)
=g-x \notin T$ and that $x \in T$ since $T$ is maximally negative semisymmetric.
Thus $M^* = T \cap (0,k)$ and $S(M^*,M,P) = T$. Clearly $P \in \Psymm(k,g-k)$
so to complete the proof it must be shown that $M \in \Asymm(k)$.
If $M_- \ne \emptyset$ then there is $x \in (0,k)$ such that $x \notin
M$ and $k-x \notin M$, and then $(g-k)+k-x = g-x \notin T$ and $x \notin
T$ which contradicts the maximality of $T$. This verifies that
$M_-$ is empty.
If $x \in M$ then $g-k+x$ is an element of $T \cap (0,g)$ larger than $g-k$
so that $g-k+x$ is not an atom of $T$ since  $g-x$ is the largest atom of $T$. 
It follows that there
is an element $y \in T$ with $y < k-x$ such that $y+g-k+x \notin T$
(note that $y$ cannot equal $k-x$ because otherwise both $k-x$
and $g-k+x$ would be elements of $T$ contradicting the assumption
that $T$ is negative semisymmetric),
and $M \in \Asymm(k)$. 
\end{proof}

By lemma~\ref{lem:LargeSmallAtom} and the theorem we have
  \begin{equation}\label{eqt:Bsigma} 
      \card{\Bsymm(g)} = \sum_{k=1}^{\floor{(g-1)/{2}}} \card{\Bsymm(g,k)} =  
      \sum_{k=1}^{\floor{(g-1)/2}} \Asigma_k\,2^{\floor{(g-2k-1)/2}} \ 
  \end{equation}
and dividing by $2^{\floor{(g-1)/2}}$ produces the next result.
	
\begin{cor}
For each $g>0$, ${\displaystyle \bsymm{g} = \sum_{k=1}^{\floor{(g-1)/2}} \Asigma_k \, 2^{-k}}$.\qed 
\end{cor}

Thus $\set{\bsymm{g}}$ is an increasing sequence, and it has a limit 
  $$\bsymm{\infty} = \sum_{k=1}^\infty \, \Asigma_k \, 2^{-k} \ .$$ 
It follows immediately that $\{\gsymm{g}\}$ is a
decreasing sequence which converges to \ $\gsymm{\infty} = 1-\bsymm{\infty}$.

\begin{cor}\label{cor:bsymmError}
For each positive integer $n$, \ $\Asigma_n \leq 3^{\floor{(n-3)/2}}$ \ and
  $$\gsymm{2n-1} - \gsymm{\infty} =  \bsymm{\infty} - \bsymm{2n-1} 
        \le \left(\frac{\sqrt{3}}{2}\right)^{n-1} .$$ 
\end{cor}

\begin{proof}
Let $M \subseteq (0,n)$ be an element of $\Asymm(n)$. Suppose that $n=2k+1$ is odd.
Then $(0,n)$ is partitioned into $k$ doubletons $\set{i, n-i}$ where $1\le i \le k$. 
Since $M_- =\emptyset$ the intersection of
$\set{i, n-i}$ with $M$ must be nonempty, and so there are three possibilities for each 
of these intersections. 
Notice also that $n-1$ cannot be an element of $M$ by condition ($\sigma$-ad2), and the
intersection of $M$ with the doubleton $\set{1,n-1}$ must be $\set{1}$.
Thus there are at most $3^{k-1}=3^{\floor{(n-3)/2}}$ possibilities for $M$.
When $n=2k$ is even, $(0,n)$ can be partitioned into $(k-1)$ doubletons
$\set{i, n-i}$ where $1\le i \le k-1$ and a singleton $\set{k}$. The intersection of
$M$ with $\set{k}$ must equal $\set{k}$ since $M_-=\emptyset$. 
By a similar argument as before there are at most $3^{k-2}=3^{\floor{(n-3)/2}}$ possibilities for $M$.
Now $\bsymm{\infty} - \bsymm{2n-1} = \sum_{k=n}^{\infty} \Asigma_k \, 2^{-k} \le
\sum_{k=n}^{\infty} 3^{\floor{(k-3)/2}}2^{-k} \le \left(\frac{\sqrt{3}}{2}\right)^{n-1}$.
\end{proof}

Some values of $\bsymm{2n-1}$ are given in table~\ref{tab:betasigma}.
Note that $\bsymm{63} = .76356\ldots$ approximates
$\bsymm{\infty}$ to within $(\sqrt{3}/2)^{31} = .0115731\ldots$
by corollary~\ref{cor:bsymmError}.
Subtracting from $1$ gives $\gsymm{63} = .23644\ldots$, which approximates 
$\gsymm{\infty}$ to within $.0115731$. Taking midpoints gives the
approximation $\gsymm{\infty} \approx .230653$ accurate to within
$.00579$.

\begin{table}[ht]
\medskip
\begin{tabular}{ccc}
\begin{tabular}{|c|c|c|c|c|}
\hline
 $n$ &$\Asp{2n-1}$& $\Asigma_n$ & $\bsymm{2n-1}$ & $R_n$ \\
\hline
   1  & 1        & 1           & 0  & 1 \\
   2  & 1        & 0           & .5  & $\infty$ \\
   3  & 2        & 1            & .5 & 1 \\
   4  & 3        & 0& .625 & $\infty$ \\
   5  & 6        & 2 & .625 & .871 \\
   6  & 10       & 0 & .6875 & $\infty$ \\
   7  & 20       & 3 & .6875 & .855 \\
   8  & 37       & 1 & .71094 & 1 \\
   9  & 73       & 7 & .71484 &.806 \\
   10 & 139      & 3 & .72852 & .896 \\
   11 & 275      & 17 & .73145 & .773\\
   12 & 533      & 7 & .73975 & .850\\
   13 & 1059     & 43 & .74146 & .749\\
   14 & 2075     & 24 & .74670 & .797\\
   15 & 4126     & 118 & .74817 & .728\\
   16 & 8134     & 74 & .75177 &.764\\
\hline
  \end{tabular}
  &&
   \begin{tabular}{|c|c|c|c|c|}
     \hline
     $n$ & $\Asp{2n-1}$ & $\Asigma_n$ & $\bsymm{2n-1}$ & $R_n$\\
     \hline
     17  & 16194& 330 & .75290 & .711\\
     18  & 32058 & 206 & .75542 & .744\\
     19  & 63910 & 888 & .75620 &  .700\\
     20  & 126932 & 612 & .75790 & .725\\
     21  & 253252 & 2571 & .75848 & .688\\
     22  & 503933 & 1810 & .75971 & .711\\
     23  & 1006056 & 7274 & .76014 & .679\\
     24  & 2004838 & 5552 & .76100 & .698\\
     25  & 4004124 & 21099 & .76134 & .671\\
     26  & 7987149 & 16334 & .76196 & .689\\
     27  & 15957964 & 61252 & .76221 & .665\\
     28  & 31854676 & 49025 & .76266 & .680\\
     29  & 63660327 & 179239 & .76285 & .659\\
     30  & 127141415 & 146048 & .76318 & .673\\
     31  & 254136782 & 523455 & .76332  & .654\\
     32  & 507750109 & 440980 & .76356 & .666\\
\hline
  \end{tabular}
\end{tabular}
\vskip10pt
\caption{Approximating $\bsymm{\infty}$, where $R_n = 1/\sqrt[n]{\Asigma_n}$.}\label{tab:betasigma}
\medskip
\end{table}

\smallskip

For a numerical monoid $M \in \Scal(g)$ let
  $$ \Gcal^\sigma(M) = \Gcal(M) \cap \Ssymm(g) = \setpres{S \in \Ssymm(g)}{$A(S) =M$} \ .$$
Note that $M$ will not be an element of $\Gcal^\sigma(M)$ unless $M$ is symmetric or
pseudosymmetric, and that $\Gcal^\sigma(M)$ may be empty.   
If $M\in \Scal(g)$ is a numerical monoid in $\Bcal(g,k)$ (which means that $g-k$ is the largest
integer in $M\cap(0,g)$) then $ \Gcal^\sigma(M) \subseteq \Bsymm(g,k)$. Therefore
\begin{displaymath}
 \card{\Gcal^\sigma(M)} \le \card{\Bsymm(g,k)} = \Asigma_k \, 2^{\floor{(g-2k-1)/2}} \le \frac{1}{3\sqrt{3}}
     \left( \frac{\sqrt{3}}{2} \right)^k \times 2^{\floor{(g-1)/2}} 
\end{displaymath} 
by corollary~\ref{cor:bsymmError}, and this 
is the symmetric analogue of the inequality in theorem~\ref{thm:AtomSet}. In particular we see that
$\card{\Gcal^\sigma(\Nbb_g)}= \card{\Gcal^\sigma(g)}$ is larger than $\card{\Gcal^\sigma(M)}$
for every numerical monoid $M \in \Scal(g)$ other than $\Nbb_g$. Taking $k=\floor{(g-1)/2}$
in the above shows that $\card{\Gcal^\sigma(\Dbb_{2k+1})} = \card{\Bsymm(g,k)} = \Asigma_k$ and that
$\card{\Gcal^\sigma(\Dbb_{2k+2})} =2\Asigma_k$, where $\Dbb_n$ is defined in 
equation~(\ref{eqt:Dbb}).

\section*{The generating function for $\set{\gsymm{k}}$}

We define $\Asymm(k)^\prime$ to be the subset of 
$\Asymm(k)$ consisting of all $\sigma$-admissible sets $M \subseteq (0,k)$ 
for which $M_+$ has at most one element.
If $M \in \Asymm(k)^\prime$ then there are two possibilities: either $k$ is odd and $M_+=\emptyset$ 
(because $M_+$ has an even number of elements whenever $k$ is odd) 
or $k$ is even and $M_+ =\set{k/2}$ (because $k/2$ must
be in $M_-$ or $M_+$ whenever $k$ is even, and $M_-=\emptyset$).
The cardinality of ${\Asymm(k)}^{\prime}$
will be denoted by $\Asp{k}$.

\begin{thm}\label{thm:Gsymm}
There is a one-to-one correspondence between the sets $\Gsymm(g)$ and ${\Asymm(g)}^{\prime}$.  
In particular, $\card{\Gsymm(g)}=\Asp{g}$ and $\gsymm{g}=\Asp{g}\,  2^{-\floor{(g-1)/2}}$.
\end{thm}

\begin{proof}
Given $M \in \Asymm(g)^\prime$ define $S=\Nbb_g \cup M^*$.
If $g$ is odd then $M_-=M_+=\emptyset$ so that $M^*=M$ and $S$ is
symmetric. If $g$ is even then $M_-=\emptyset$, $M_+ = \set{g/2}$ and $M^*=M-\set{g/2}$
which implies that $S$ is pseudosymmetric. In either case, $S$ is an element of $\Ssymm(g)$
and we can write $M^* = M - \set{g/2}$. Suppose $s$ is a large small atom in $S$. Then
$s \in M^* \subseteq M$ and by ($\sigma$-ad2) there exists $y \in M^*$ such that $s+y < g$
and $s+y \notin M$. Note that $s+y \notin M-\set{g/2} = M^*$ since $s$ is large. Thus
$s+y \notin S$ in contradiction of the assumption that $s \in A(S)$. This shows that
$M \mapsto S$ is an injective function from $\Asymm(g)^\prime$ to $\Gsymm(g)$. If $S \in \Gsymm(g)$
then it is not hard to check that $S=\Nbb_g \cup M^*$ where $M = (S \cap (0,g))^* \in 
\Asymm(g)^\prime$, completing the proof.
\end{proof}

\begin{thm}\label{thm:symmetricrecursion}
For each $k \geq 1$, $\Asp{2k} = \Asp{2k-1}$ and  $\Asp{2k+1} = 2\Asp{2k} - \Asigma_{k}$.
\end{thm}

\begin{proof}
The first statement follows from lemma~\ref{lem:EvenOddSymm} and theorem~\ref{thm:Gsymm}. 
For the second we have
  $$\Asp{2k+1} = \card{\Gsymm(2k+1)} =  \card{\Ssymm(2k+1)}- \card{\Bsymm(2k+1)} = 
       2^k - \sum_{\ell=1}^k \, \Asigma_\ell 2^{k-\ell}$$
using equation~(\ref{eqt:Bsigma}), and
$\Asp{2k} =  2^{k-1} - \sum_{\ell=1}^{k-1} \, \Asigma_\ell 2^{k-\ell-1}$
by a similar computation. Combining these two equations gives $2\Asp{2k} - \Asp{2k-1} = \Asigma_k$.
\end{proof}

\begin{figure}[hbt]
\bigskip
    \centering
    \includegraphics[width=3.5in]{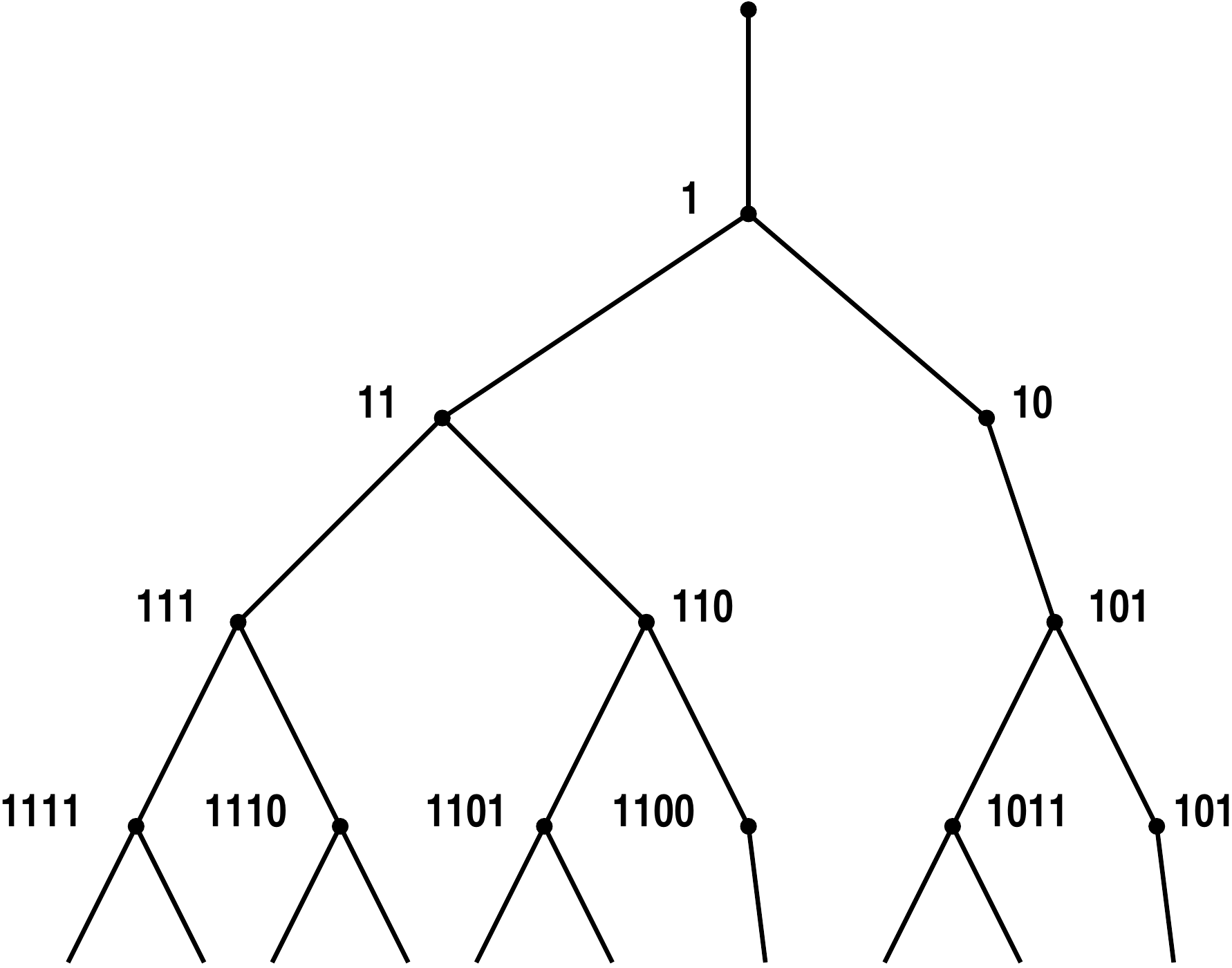}
    \caption{\textsl{Rooted tree for \ $\bigcup \setpres{\Gsymm(2k+1)}{$k \in \Nbb$}$.}}\label{fig:sigmaRootedTree}
\medskip
\end{figure}

As in figure~\ref{fig:RootedTree} we may view $\bigcup\setpres{\Gsymm(2k+1)}{$k \in \Nbb$}$
as the vertices of a downward opening rooted tree.
Here each element of $\Gsymm(2k-1)$ (represented by a vertex in the $k$th level of the tree)
will spawn either one or two elements of $\Gsymm(2k+1)$ (represented by vertices
at the $(k+1)$st level).
If $S \in \Gsymm(2k-1)$ and $Q$ is either $\set{k}$ or $\set{k+1}$
then let $S(Q)$  be defined by equation~(\ref{eqt:spawn}) and observe that $S(Q) \in \Ssymm(2k+1)$.
As before $S(Q)$ has no small atoms when $Q= \set{k}$ but 
$k+1$ may be a small atom for $S(Q)$ when $Q = \set{k+1}$.
From this we see that $\card{\Gsymm(2k+1)} = \Asp{2k+1}$ is two times $\card{\Gsymm(2k-1)}=\Asp{2k-1}$
minus the number of elements of $\Gsymm(2k-1)$ which spawn only one element of $\Gsymm(2k+1)$, and
it follows from theorem~\ref{thm:symmetricrecursion} that
the number of elements of  $\Gsymm(2k-1)$ which spawn only one element of $\Gsymm(2k+1)$ equals $\Asigma_k$.      
Figure~\ref{fig:sigmaRootedTree} shows the first few levels of the rooted tree. 
In this illustration 
a labeling sequence $\alpha = (\alpha_1, \ldots , \alpha_k)$ with entries in $\Zbb_2$ 
represents the same numerical set in $\Gsymm(2k+1)$as the 
$2 \times k$ matrix 
$\left( \begin{smallmatrix} 
   \alpha_1 &\cdots &\alpha_k \\ \alpha_{1}^*  &\cdots &\alpha_{k}^* 
\end{smallmatrix} \right)$ 
represented in figure~\ref{fig:RootedTree}, where $\alpha_i^* = 1 -\alpha_i$.  
Call the sequence $\alpha$ {\it bivalent} if
there is an integer $\ell$ with $1 \le \ell \le k$ such that $\alpha_\ell
=\alpha_{k+1-\ell}=1$. Then $\alpha$ represents an element of $\Gsymm(2k+1)$
if and only if $(\alpha_1, \ldots, \alpha_{i-1})$ is bivalent whenever $\alpha_i=0$.

If we associate $\alpha = (\alpha_1, \ldots , \alpha_k)$ with the finite set $F(\alpha)=\set{0}\cup\setpres{i}{$\alpha_i=1$}$
then $\alpha$ is bivalent if and only if $k+1 \in F(\alpha)+F(\alpha)$. Furthermore,
$\alpha$ represents an element of $\Gsymm(2k+1)$ precisely when $F(\alpha)$ is an
`additive $2$-basis for $k$' (which means that $[0,k) \subseteq F(\alpha)+F(\alpha)$). 
Thus $\Asp{2k+1}$ equals the number of subsets of $[0,k)$ which are additive $2$-bases for $k$. 
By theorem~\ref{thm:Gsymm}, $\Asp{2k+1}=\gsymm{k} 2^k$ which
is asymptotically equal to $\gsymm{\infty} 2^k$ where $\gsymm{\infty} \approx .230653$ as described above.
The study of finite additive $2$-bases for $k$ has a long history, especially in relation to the 
determination of bounds for the smallest cardinality of such bases.  
The introduction of  \cite{GN} has a nice overview of this. 
The first $19$ terms of the sequence $\{ \Asp{2k-1} \}$  have been posted at  
\cite{S} by M.~Torelli as sequence number A008929. The paper \cite{T} describes some 
related sequences. (In that paper a finite additive 2-basis is called a 
(finite) `Goldbach sequence'.)

\medskip

Let $g^\sigma(z)$ and $f^\sigma(z)$ be the analytic functions defined by 
  $$ g^\sigma(z) = \sum_{k=1}^\infty \, \Asigma_k z^k \quad \text{and} \quad  
     f^\sigma(z) = \sum_{k=1}^\infty \, \Asp{k} z^k \ .$$
  
\begin{cor}\label{cor:fgsigma}
The functions $f^\sigma(z)$ and $g^\sigma(z)$ satisfy the relation
$$\left(2z^2-1\right)f^\sigma(z) = z(z+1)\left(g^\sigma(z^2)-1\right) \ .$$
\end{cor}

\begin{proof}
First observe that
  \begin{equation}\label{eqt:fsigma}
  f^\sigma(z) = \sum_{k=1}^\infty \, \Asp{2k-1}z^{2k-1} +  \sum_{k=1}^\infty \, \Asp{2k}z^{2k}
    = (z+1) \sum_{k=1}^\infty \, \Asp{2k-1}z^{2k-1}
  \end{equation}
by theorem~\ref{thm:symmetricrecursion}. From the same theorem, $\Asigma_k = 2\Asp{2k-1} -  \Asp{2k+1}$
and
 $$\begin{aligned}
    g^\sigma(z^2) &= \sum_{k=1}^\infty \, 2\Asp{2k-1}z^{2k} -  \sum_{k=1}^\infty \, \Asp{2k+1}z^{2k}\cr
   &= 2z\sum_{k=1}^\infty \Asp{2k-1}z^{2k-1} - \frac{1}{z}\sum_{k=1}^\infty\Asp{2k-1}z^{2k-1} + \Asp{1}\cr
   &= \left( 2z - \frac{1}{z}\right) \sum_{k=1}^\infty\Asp{2k-1}z^{2k-1} + 1 = \frac{2z^2-1}{z(z+1) }f^\sigma(z) + 1 . 
   \end{aligned}$$
\end{proof}

\begin{cor} The analytic function $f^\sigma(z)$ has singularities at 
$z=\pm 1/\sqrt{2}$ and its radius of convergence at the origin equals
$1/\sqrt{2}$. Except for $z=\pm 1/\sqrt{2}$
and possibly for $z=-1$, the singularities
of $f^\sigma(z)$ coincide with those of $g^\sigma(z^2)$ and 
$f^\sigma(z)(2z^2-1)$ is continuous on the closed disk $\abs{z} \le 1/\sqrt{2}$.
\end{cor}

\begin{proof}
Since $g^\sigma(1/2) = \bsymm{\infty}$ the power series
$\sum_{k=1}^\infty \, \Asigma_k z^k$ converges absolutely 
and $g^\sigma(z)$ is continuous on $\abs{z} \le 1/2$.
By corollary~\ref{cor:fgsigma} 
$f^\sigma(z)$ has singularities 
at $z=\pm 1/\sqrt{2}$ 
and $f^\sigma(z)$ has radius of convergence $1/\sqrt{2}$ at the origin.
The last property also follows immediately from corollary~\ref{cor:fgsigma}. 
\end{proof}

Since $0 <  \Asp{k} < \Asigma_k$, the series $\sum_{k=1}^\infty \, \Asigma_k \left(1/\sqrt{2}\right)^k$ 
diverges by comparison with $\sum_{k=1}^\infty \, \Asp{k} \left( 1/\sqrt{2}\right)^k$, and so
the radius of convergence of $g^\sigma(z)$ at the origin must be
between $1/2$ and $1/\sqrt{2}$. 
The root test would equate this radius of convergence with the limit infimum of 
$R_n = 1/\sqrt[n]{\Asigma_n}$. This value seems to be larger than $1/2$ by
the data in the last column of table~\ref{tab:betasigma},
but we have not been able to ascertain this.

Let $h^\sigma(z) = \sum_{k=1}^\infty \, \gamma^\sigma_k z^k $ be the generating function for 
$\set{\gamma^\sigma_k}$. 
Using equation~(\ref{eqt:fsigma}) and corollary~\ref{cor:fgsigma} we have
 $$\begin{aligned}
  h^\sigma(z) &= \sum_{k=1}^\infty \frac{\Asp{k}}{2^{\floor{(k-1)/2}}} \, z^k =
         \sum_{k=1}^\infty \frac{\Asp{2k-1}}{2^{k-1}}\, z^{2k-1} + \sum_{k=1}^\infty 
         \frac{\Asp{2k}}{2^{k-1}}\, z^{2k}  \cr
  &= \sqrt{2}\Big( z+1 \Big) \sum_{k=1}^\infty \Asp{2k-1} \left(\frac{z}{\sqrt{2}}\right)^{2k-1} 
  = 2\left( {z+1}/{z+\sqrt{2}}\right)  \,f^\sigma\left({z}/{\sqrt{2}}\right) \cr
  &= \left( \frac{z}{z-1} \right) \left( g^\sigma\left({z^2}/{2}\right) -1\right) . 
  \end{aligned}$$
Therefore $h^\sigma(z)$ has radius of convergence $1$ at the origin and has a singularity at
$z=1$.
If the radius of convergence of $g^\sigma(z)$ at the origin is
larger than $1/2$ then $z=1$ is the only singularity of $h^\sigma(z)$ inside a circle with radius larger
than $1$ centered at the origin and 
$h^\sigma(z)$ would have a simple pole at $z=1$ with residue
  $$\lim_{z\to 1}(z-1)h^\sigma(z) = \lim_{z \to 1} z(g^\sigma(z^2/2)-1) = g^\sigma(1/2) -1 = -\gsymm{\infty}  
     \ .$$

\vfil


\begin{thebibliography}{999}
\bibitem[A]{A}
   J.L.Alfonsin,
   \textit{The Diophantine Frobenius Problem},
   \textsl{Oxford University Press}, 2005.
\bibitem[AM]{AM}
   E.~Antokoletz and A.~Miller,
   \textit{Symmetry and factorization of numerical sets and monoids},
   \textsl{J.~of Algebra} 247 (2002), 636--671.
\bibitem[B]{B}
   J.~Backelin,
   \textit{On the number of semigroups of natural numbers},
   \textsl{Mathematica Scandinavica} 66 (1990), 197--215.
\bibitem[BF]{BF}
   V.~Barrucci and R.~Fr\"oberg,
   \textit{One dimensional almost Gorenstein rings},
   \textsl{J.~of Algebra} 188 (1997), 418--442.
\bibitem[FGH]{FGH}
   R.~Fr\"oberg, C.~Gottlieb, and R.~H\"aggkvist
   \textit{On numerical semigroups},
   \textsl{Semigroup Forum} 35 (1987), 63--83.
\bibitem[GN]{GN}
   S.~G\"unt\"urk and M.~Nathanson,
   \textit{A new upper bound for finite additive bases},
   \textsl{Acta Arith.} 124 (2006), 235--255.
\bibitem[S]{S}
    N. J. A. Sloane, \textit{The On-Line Encyclopedia of Integer Sequences}, published electronically at \textsl{http://www.research.att.com/}\verb+~+\textsl{njas/sequences}, (2008).
\bibitem[T]{T}
   M.~Torelli,
   \textit{Increasing integer sequences and Goldbach's conjecture},
   \textsl{Inf.~Theor.~Appl.} 40 (2006), 107--121.
\end{thebibliography}
\end{document}